\documentclass[11pt]{article}
\usepackage{color}
\usepackage{amsbsy}
\usepackage{amsmath}
\usepackage{amssymb}
\usepackage{latexsym}
\usepackage{comment}
\usepackage{graphicx}
\usepackage{amsfonts}
\usepackage{url}
\usepackage{hyperref}
\usepackage[hyphenbreaks]{breakurl}
\usepackage{enumitem}
\usepackage{bbm}

\usepackage[authoryear]{natbib}

\usepackage{geometry}

\def\singlespace{\def\baselinestretch{1}\@normalsize}

\renewcommand{\baselinestretch} {1.2}
\makeatletter 
\def\singlespace{\def\baselinestretch{1}\@normalsize}

\@addtoreset{equation}{section}

\newtheorem{theorem}{Theorem}

\newtheorem{lemma}{Lemma}
\newtheorem{remark}{Remark}
\newtheorem{assumption}{Assumption}

\newtheorem{definition}{Definition}

\newcommand{\norm}{\vert\kern-0.3ex\vert\kern-0.3ex\vert} 
\newcommand{\Bignorm}{\Big\vert\kern-0.4ex\Big\vert\kern-0.4ex\Big\vert} 

\def\tilde{\widetilde}
\def\hat{\widehat}
\def\nullH0{\mathrm{null}}

\begin{document}

\title{On Controlling the False Discovery Rate in Multiple Testing of the Means of Correlated Normals Against Two-Sided Alternatives}

\author{Sanat K. Sarkar\thanks{Department of Statistics, Operations, and Data Science, Temple University, USA. Email: sanat@temple.edu. The research is supported by NSF grant DMS 2210687}}
\maketitle

\begin{abstract} This paper revisits the following open question in simultaneous testing of multivariate normal means against two-sided alternatives: Can the method of Benjamini and Hochberg (BH, 1995) control the false discovery rate (FDR) without imposing any dependence structure on the correlations? The answer to this question is generally believed to be yes, and is conjectured so in the literature since results of numerical studies investigating the question and reported in numerous papers strongly support it. No theoretical justification of this answer has yet been put forward in the literature, as far as we know. In this paper, we offer a partial proof of this conjecture. More specifically, we  consider the following two settings - (i) the covariance matrix is known and
(ii) the covariance matrix is an unknown scalar multiple of a known matrix
 - and prove that in each of these settings a BH-type stepup method based on some weighted versions of the original $z$- or $t$-test statistics controls the FDR.

\bigskip
\end{abstract}


\section{Introduction} The false discovery rate (FDR) introduced by \cite{Benjamini1995} is a powerful notion of an overall measure of type I error in multiple testing. With multiple testing being an ubiquitous inferential tool in statistical investigations arising in modern scientific research, the method of \cite{Benjamini1995} designed to control this error rate, popularly known as the BH method, is now one of the most commonly used multiple testing methods. Despite the immense popularity of the BH method, the use of it as a valid FDR controlling method, unfortunately, is often questioned in many practical applications, since conditions ensuring such validity are not often met in the multiple testing scenarios encountered in those applications. One such scenario involves multiple testing of the means of correlated normal  random variables with a positive definite correlation matrix. The BH method provably controls the FDR in this scenario when the alternatives are all one-sided, as long as the correlations are non-negative [see, e.g., \cite{Benjamini2001}, \cite{Blanchard2008},
\cite{Finner2007}, \cite{Sarkar2002}]. However, no such result is known in the literature about its FDR control when the alternatives are all two-sided. Testing against two-sided alternatives is most often scientifically more meaningful than testing against one-sided alternatives. Thus, the applicability of the BH method in many modern statistical applications is limited without being theoretically verified as a valid FDR controlling method in the aforementioned two-sided testing scenario. This has been one of the motivations that led to the recent upsurge of research bypassing the use of the BH method and developing alternative methods [\cite{Barber2015}, \cite{Fithian2022}, \cite{Sarkar2022}] that can provably control the FDR under the same multiple testing scenarios.

Of course, there is a strong belief among multiple testing researchers, and it is often conjectured, that the BH method can indeed control the FDR for multiple testing of correlated normal means against two-sided alternatives no matter what the correlation matrix is. As \cite{Benjamini2010} remarked:
\begin{quote}
{\it The modification to general dependence is often not needed: convincing simutheoretical evidence indicates that the same holds for two-sided $z$-tests with any correlation structure [\cite{Reiner-Benaim2007}], but the theory awaits a complete proof}.
\end{quote}
A strong support for this belief can also be seen from results of numerical studies carried out in many papers that investigated the BH method's FDR control for such testing problem. So, in light of the aforementioned recent developments of alternative methods, proving this conjecture and thus making the BH method a theoretically valid FDR controlling method as a relevant competitor for these other methods seems an urgent and important undertaking.

This paper presents a partial proof of this conjecture with a known but arbitrary correlation matrix in terms of a weighted version of the BH method. More specifically, we prove the conjecture considering a BH-type stepup method in which the $z^2$-values in two-sided $z$-tests [mentioned in the above remark of \cite{Benjamini2010}] are weighted according to the extent to which each $z$ correlates with the others. Expanding the arguments used in this proof, we then give our proof for the BH-type stepup method involving similarly weighted $t^2$-values in two-sided $t$-tests.
These proofs, as special cases, establish the validity of the Simes global test involving the aforementioned weighted $z$- and $t$-test statistics. 
See \cite{Finner2017} for the conjecture made for the Simes global test involving the un-weighted $z$- and $t$-test statistics.

The paper is organized as follows. Section 2 presents some basic results in terms of formulas for the FDR of a stepup test and conditions under which the FDR can be controlled. Our proposed methods and the main results associated with  them  are given in Section 3. In Section 4, the methods are expounded in the contexts of equi-correlated multivariate normal and variable selection under linear regression model. A novel BH-type FDR controlling procedure for variable selection is produced. The paper concludes with some additional remarks in Section 5.

\section{Preliminaries} In this section, we present some basic notations, formulas, and assumptions associated  with the BH or a closely related stepup  method, before having further discussions setting the stage for our main results in the next section.

Given a set of $d$ null hypothesis $H_i, \; i=1, \ldots, d$, to be tested simultaneously using their respective $p$-values or some increasing functions of them, $P_i$, $i=1, \ldots, d$, the BH method is a stepup test applied to the $P_i$'s with critical constants $\alpha_i=  i \alpha_1$, $i=1, \ldots, d$; that is, it finds $R = \max_{1 \le i \le d} \{i: P_{(i)} \le \alpha_i \}$, and rejects $H_i$ for all $i$ such that $P_i \le P_{(R)}$, provided the maximum exists; otherwise, it rejects none. It is designed to control
\begin{eqnarray} {\rm FDR} = \mathbb{E} (\textrm{FDP}), \; \mbox{where FDP} \; \mbox{(False Discovery Proportion)}\; = \frac{V}{R \vee 1}, \end{eqnarray} ($a \vee b = \max(a,b)$), with $V$ and $R$ being the numbers of falsely rejected and rejected null hypotheses, respectively, at level $d\alpha_1$.

The following lemma provides an explicit expression for the FDR of a stepup test with any set of critical constants:

\begin {lemma} Let $I_0= \{j: H_j\; \mbox{is true} \}$. Then, the FDR of a stepup test with critical constants $0< \alpha_1 \le \cdots \le \alpha_{d}< 1$ applied to the $P_i$'s is given by
\begin {eqnarray}\label{FDR1}
\textrm {FDR} & = & \sum_{i \in I_0}\mathbb{E} \left \{ \frac {\mathbbm{1}(P_i \le \alpha_{R_{-i}+1})}{R_{-i}+1} \right \} \nonumber \\ & = & \sum_{i \in I_0} \sum_{r=0}^{d-1} \mathbb{E} \left \{ \frac{\mathbbm{Pr} \left (P_i \le \alpha_{r+1}~|~\boldsymbol{P}_{-i} \right ) \mathbbm{1} \left (R(\boldsymbol{P}_{-i}) = r \right )}{r+1} \right \}, \end{eqnarray}
where $\boldsymbol{P}_{-i} = (P_1, \ldots, P_d)\setminus\{P_i\}$ and $ R_{-i}\equiv R(\boldsymbol{P}_{-i}) = \max_{1\le j \le d-1}\{j:P_{(j)\setminus\{i\}} \le \alpha_{j+1} \},$ with $P_{(1)\setminus\{i\}} \le \cdots \le P_{(d-1)\setminus\{i\}}$ being the ordered components of $\boldsymbol{P}_{-i}$. \end{lemma}

See Sarkar (2002, 2008) for this formula, although similar formulas do appear in other papers as well [e.g., \cite{Benjamini2001}, \cite{Blanchard2008}, \cite{Finner2007}].

When the $P_i$'s are independent, it is immediate from this formula that the FDR of the BH method based on these $P_i$'s equals $\sum_{i \in I_0}\mathbbm{Pr} (P_i \le \alpha_1)$, and hence is controlled at $|I_0|\alpha_1 \le d \alpha_1$, where $|I_0|$ is the cardinality of $I_0$, under the following assumption:

\begin{assumption} For each $i \in I_0$, $P_i$ is stochastically larger than the random variable with $U(0,1)$ distribution. \end{assumption}

When the $P_i$'s are not independent, an approach to finding the condition under which the BH method or a closely related stepup method can continue to control the FDR is to re-write the above formula in an alternative form that can reveal the type of dependence for the $P_i$'s one would need to prove the FDR control. One such formula is
\begin {eqnarray}\label{FDR2}
\textrm {FDR} & = & \sum_{i \in I_0}\mathbbm{Pr} (P_i \le \alpha_1) + \nonumber \\ & & \sum_{i \in I_0} \sum_{r=1}^{d-1} \mathbb{E} \left \{ \mathbbm{Pr} \left (R(\boldsymbol{P}_{-i}) \ge r ~|~P_i \right ) \left [ \frac{\mathbbm{1}(P_i \le \alpha_{r+1})}{r+1}- \frac{\mathbbm{1}(P_i \le \alpha_{r})}{r} \right ] \right \},
\nonumber \\ \end{eqnarray}
(\cite{Sarkar2002}). The set $\left \{\boldsymbol{P}_{-i}: R(\boldsymbol{P}_{-i}) \ge r \right \}$ is decreasing in $\boldsymbol{P}_{-i}$, for any fixed $r=1, \ldots, d-1$. This is the crux of a proof of the FDR control under Assumption 1, and led \cite{Benjamini2001} and others to consider making the following assumption on the dependence structure of the $P_i$'s ensuring the above decreasing property:

\begin{assumption} The $P_i$'s are positively regression dependent on the subset (PRDS) of  $P_i$'s corresponding to the null hypotheses, i.e., \begin {eqnarray}
\mathbb{E} \left \{ \phi(\boldsymbol{P}_{-i})~|~ P_{i} \right \} \; \uparrow \; (\mbox{or} \; \downarrow ), \; P_{i} \; \mbox{for each} \; i \in I_0, \end {eqnarray}
and for any co-ordinatewise increasing (or decreasing) function of $\boldsymbol{P}_{-i}$. \end{assumption}

The PRDS is a positive dependence condition that, being satisfied in many practical scenarios including the ones where $p$-values are generated from normal test statistics with non-negative correlations, has now been accepted as the only positive dependence condition under which the BH method can control the FDR in a non-asymptotic setting. Unfortunately, however, it does not capture the positive dependence structure exhibited by test statistics, and hence by the corresponding $p$-values, arising in many other and relatively more important scenarios.  Among them are those where the test statistics have folded multivariate normal or folded multivariate $t$ distribution. These statistics arise in the context of multiple testing of the means of correlated normals with known or with unknown variances against two-sided alternatives. The tools or distributional properties of multivariate normal for checking the PRDS condition in the case of one-sided testing problems no longer work for folded multivariate normal or $t$ test statistics used for two-sided testing problems. For instance, when the test statistics are multivariate normal, the conditional distribution of $\boldsymbol{P}_{-i}$ given $P_i$ stochastically increases with $P_i$ when the correlations are nonnegative, from which the PRDS condition can be verified for one-sided testing problems. Similarly, when the partial correlations are all positive, multivariate normal is known to be totally positive of order two (MTP$_2$) that implies the PRDS condition; see \cite{Karlin1980} and \cite{Sarkar2002}. These tools and distributional properties don't work for the BH method based on folded multivariate normal or folded $t$-test statistics.

Thus, whether or not the BH method involving $z^2$- or $t^2$-values, or some increasing functions of them, provably controls the FDR for multiple testing of the means of multivariate normal
with arbitrary but known positive definite correlation matrix in a non-asymptotic setting has remained one of the important open problems in multiple testing.  The next section resolves this problem in terms of what we call weighted BH methods.

\begin{remark} \rm It is important to clarify at this point what we mean by a weighted BH method when the weighting scheme, given some weights $w_i>0$, $i=1, \ldots, d$, involves the underlying test statistics, not the $p$-values as typically assumed in the literature. Let $T_i$ be the test statistic generating the $p$-value $P_i = \bar{F}_0(T_i)$, for $i=1, \ldots, d$, using its null survival function $\bar{F}_0$, and $w_i^{-1}T_i$, for $i=1, \ldots, d$, be  weighted versions of the $T_i$'s. Then, we refer to the BH method applied to the $\tilde{P}_i = \bar{F}_0(w_i^{-1}\bar{F}_0^{-1}(P_i))$, $i=1, \ldots, d$, with critical constants $i\alpha_1$, for $i=1, \ldots  d$, as a weighted BH method based on the $P_i$'s. Alternatively, if we focus on the BH method in terms of the test statistics, rather than the $p$-values; that is, find $R= \min\{i: T_{(i)} \ge \bar{F}_0^{-1} ((d-i+1)\alpha_1) \}$, having ordered the $T_i$'s as $T_{(1)} \le \cdots \le T_{(d)}$, and reject $H_i$ for $i$ such that $T_i \ge T_{(R)}$, provided the minimum exists; otherwise, rejects none, then we refer to this method with the $T_i$'s replaced by their weighted versions, defined using some weights assigned to them, as a weighted BH method.
\end{remark}

\section{Main Results} \label{s3}

We present in this section our main results of this paper, Theorems 1 and 2, related to the aforementioned two-sided tests for multiple testing of the means of multivariate normal. In these theorems, the correlation matrix is assumed known without exhibiting any specific dependence structure. Theorem 1 presents our proposed FDR controlling BH-type stepup method based on weighted $z^2$-values and Theorem 2 presents the same based on weighted $t^2$-values.

\subsection{Weighted $z^2$ values} Suppose that we have a $d$-dimensional random vector  $\boldsymbol{X}= (X_1, \ldots, X_d)^{\prime} \sim \mathcal{N}_d(\boldsymbol{\mu}, \boldsymbol{\Sigma})$, with an unknown mean vector $\boldsymbol{\mu} = (\mu_1, \ldots, \mu_d)^{\prime}$ and a known positive definite covariance matrix $\boldsymbol{\Sigma} = ((\sigma_{ij}))$, and that our problem is to test $H_i: \mu_i=0$ against $\mu_i \neq 0$, simultaneously for $i=1, \ldots, d$, subject to a control of the FDR at $\alpha$.

Let $\bar{\Psi}_{n} = 1 -{\Psi}_{n}$, where $\Psi_{n}$ denotes the cdf of $\chi_{n}^2$, the central chi-squared random variable with $n$ degrees of  freedom. Then, $P_i = \bar{\Psi}_{1}(Z_i^2)$, $i=1, \ldots, d$, where $Z_i= X_i/\sqrt{\sigma_{ii}}$, are the $p$-values in the two-sided $z$-tests. Instead of applying the BH method to these $p$-values, we consider applying it to the following increasing functions of them using the critical constants $i\alpha_1$, $i=1, \ldots, d$, with some appropriately chosen $\alpha_1$ depending on the level $\alpha$ at which the FDR is to be controlled:
\begin{eqnarray}\label{p-val1} \tilde{P}_i = \bar{\Psi}_{1} \left ( w_i^{-1}\bar {\Psi}_{1}^{-1} (P_i) \right ), \; \mbox{where} \; w_i = 1-R_i^2, \; i=1, \ldots, d, \end{eqnarray} and $R_i^2$ is the squared multiple correlation between $Z_i$ and $(Z_1, \ldots, Z_d)\setminus\{Z_i\}$. In other words, we consider assigning to each $Z_i^2$ the weight $w_i$ that reflects the extent to which $Z_i$ is correlated with the others, and having defined $Y_i = w_i^{-\frac{1}{2}}Z_i$, for $i=1, \ldots, d$, we  run the BH on the $Y_i^2$'s, the weighted $z^2$-values, using the critical constants $\bar{\Psi}_1^{-1}((d-i+1)\alpha_1)$, $i=1, \ldots, d$ (as explained in Remark 1). This is our proposed weighted BH method involving two-sided $z$-tests. The fact that it controls the FDR is stated in the following theorem.

\begin{theorem} \label{theorem1} The FDR of the BH method applied to the $\tilde{P}_i$'s in (\ref{p-val1}) using the critical constants $i\alpha_1$, $i=1, \ldots, d$, with $\alpha_1$ satisfying
\begin {eqnarray} \label{cond_theorem1}
& & \sum_{i=1}^d \bar{\Psi}_{1} \left ( w_i \bar {\Psi}_{1}^{-1} (\alpha_1) \right ) = \alpha, \end{eqnarray} is controlled at $\alpha$. \end{theorem}

Before proving this theorem, let us present two lemmas that will play key roles in our proof of the theorem.

The first lemma (Lemma 2) presents a useful result related to the distribution of $\boldsymbol{Y} = (Y_1, \ldots, Y_d)^{\prime}$. While deriving this distribution, we first note that $\boldsymbol{Z} = (Z_1, \ldots, Z_d)^{\prime} \sim N_{d} (\boldsymbol{\nu}, \boldsymbol{\Sigma})$, where $\boldsymbol{\nu} = (\nu_1, \ldots, \nu_d)$, with  $\nu_i=\mu_i/\sqrt{\sigma_{ii}}$, for $i=1, \ldots, d$, and that $\boldsymbol{\Sigma}$ can be assumed to be the correlation matrix. Also, the $i$th diagonal entry of $\boldsymbol{\Sigma}^{-1}$ equals $(1-R_i^2)^{-1}=w_i^{-1}$. Hence, $\boldsymbol{Y} = \textrm{diag} \{w_i^{-\frac{1}{2}} \}\boldsymbol{Z} \sim N_{d}(\boldsymbol{\delta}, \boldsymbol{\Gamma}^{-1})$, where $\boldsymbol{\delta} = (\delta_1, \ldots, \delta_d)^{\prime}$, with $\delta_i= \nu_{i}/\sqrt{w_i}= \mu_i/\sqrt{\sigma_{ii}(1-R_i^2)}$, for $i=1, \ldots, d$, and $\boldsymbol{\Gamma} = \textrm{diag}\{ w_i^{\frac{1}{2}}  \}\boldsymbol{\Sigma}^{-1} \textrm{diag} \{ w_i^{\frac{1}{2}}\},$ whose $i$th diagonal entry is $1$.

In terms of the following notations:
\begin{itemize}
\item [$\bullet$] $\boldsymbol{Y}_{-i}$: the $(d-1)$-dimensional sub-vector of $\boldsymbol{Y}$ without its $i$th entry ${Y}_i$,\
\item [$\bullet$] $\boldsymbol{\delta}_{-i}$: the $(d-1)$-dimensional sub-vector of $\boldsymbol{\delta}$ without its $i$th entry ${\delta}_i$,\
\item [$\bullet$] $\boldsymbol{\gamma}_{-i,i}$: the $(d-1)$-dimensional sub-vector of the $i$th column of $\boldsymbol{\Gamma}$ without its $i$th entry $1$, and \
\item [$\bullet$] $\boldsymbol{\Gamma}_{-i,-i}$: the $(d-1) \times (d-1)$ principal sub-matrix of $\boldsymbol{\Gamma}$ without its $i$th row and $i$th column,
\end{itemize}
we then see that
\begin {eqnarray} Y_i - \delta_i~|~\boldsymbol{Y}_{-i} & \sim & N (- \boldsymbol{\gamma}_{-i,i}^{\prime} (\boldsymbol{Y}_{-i}- \boldsymbol{\delta}_{-i}), 1) \nonumber \\ \boldsymbol{Y}_{-i} & \sim & N_{d-1}\left (\boldsymbol{\delta}_{-i}, [ \boldsymbol{\Gamma}_{-i,-i} -\boldsymbol{\gamma}_{-i,i}\boldsymbol{\gamma}_{-i,i}^{\prime}]^{-1}\right ), \nonumber \end{eqnarray}
yielding the following lemma.

\begin{lemma} When $\mu_i=0$, \begin {eqnarray} Y_i^2~|~\boldsymbol{Y}_{-i} \sim {\chi^{\prime}_{1}}^2 (\lambda_i(\boldsymbol{Y}_{-i})), \end{eqnarray} non-central chi-square with 1 degree of freedom and the non-centrality parameter $\lambda_i (\boldsymbol{Y}_{-i}) = (\boldsymbol{\gamma}_{-i,i}^{\prime} (\boldsymbol{Y}_{-i}- \boldsymbol{\delta}_{-i}))^2$.
\end{lemma}

The next lemma presents a useful result related to non-central chi-square distribution.

\vskip 6pt

\begin{lemma}  \begin {eqnarray} \frac{\mathbbm{Pr}[{\chi^{\prime}}_{n}^2 (\lambda) \ge  \bar{\Psi}^{-1}_{n}(u)]} {u} \downarrow u\in (0,1), \nonumber \end{eqnarray}  for any fixed $n, \lambda > 0$. \end{lemma}

{\it Proof}. Using the fact that ${\chi^{\prime}_{n}}^2 (\lambda) \stackrel{d} = \mathbb{E}_{J}(\chi_{n+2J}^2)$, where the expectation is taken with respect to $J \sim {\rm Poisson}\left (\frac{\lambda}{2} \right )$, we first note that \begin{eqnarray} \label{Non-Cent}\mathbbm{Pr}[{\chi^{\prime}_{n}}^2 (\lambda) \ge  {\bar{\Psi}^{-1}}_{n}(u)] = \mathbb{E}_{J} \left [ \bar{\Psi}_{n+2J}(\bar{\Psi}_n^{-1}(u) \right ]. \end{eqnarray}

Considering the function $g(u) = {\bar{\Psi}}_{n+h}({\bar{\Psi}}^{-1}_{n}(u)), \; u \in (0,1)$, for any fixed $n, h >0$, and noting that the density $\psi_n(y)$ of $\Psi_n(y)$ is $\propto e^{-\frac{1}{2}y} y^{\frac{n}{2}-1} \mathbbm{1}(y >0); n >0$, we then see that
$$ \frac{dg(u)}{du} = \frac{\psi_{n+h}(\bar{\Psi}^{-1}_{n}(u))}
{\psi_{n}(\bar{\Psi}^{-1}_{n}(u))}
\propto \left(\bar{\Psi}^{-1}_{n}(u)\right )^{\frac{h}{2}} \downarrow u \in (0,1);$$ that is, $g(u)$ is concave in $u\in(0,1)$. Since $g(0)=0$, the concavity of $g$ implies that $\frac{g(u)}{u} \downarrow u \in (0,1)$. Using this in (\ref{Non-Cent}), we have the proof of Lemma 3.

\vskip 6pt

{\it Proof of Theorem 1.} We use the formula in (\ref{FDR1})  with $\tilde{P}_i$ and $\tilde{\boldsymbol{P}}_{-i}$ being written in terms $Y_i$ and $\boldsymbol{Y}_{-i}$, respectively. From Lemma 2, we then see that the FDR of a stepup test applied to the $\tilde{P}_i$'s in (\ref{p-val1}) with any set of critical constants $\alpha_i, i=1, \ldots, d$, is given by

\begin {eqnarray}\textrm{FDR} & = & \sum_{i \in I_0} \sum_{r=0}^{d-1} \mathbb{E} \left \{ \frac{\mathbbm{Pr} [Y_i^2 \ge \bar{\Psi}_{1}^{-1}(\alpha_{r+1})~|~ \boldsymbol{Y}_{-i}]\mathbbm{1} \left (R(\boldsymbol{Y}_{-i}) = r \right )}{r+1} \right \} \nonumber \\
& = & \sum_{i \in I_0} \sum_{r=0}^{d-1} \mathbb{E} \left \{ \frac{\mathbbm{Pr}[{\chi^{\prime}_{1}}^2 (\lambda_i(\boldsymbol{Y}_{-i})) \ge  \bar{\Psi}_{1}^{-1}(\alpha_{r+1})]}{r+1} \mathbbm{1} \left (R(\boldsymbol{Y}_{-i}) = r \right ) \right \}\nonumber \\
& \le & \sum_{i \in I_0} \sum_{r=0}^{d-1} \mathbb{E} \left \{ \frac{\alpha_{r+1}}{r+1} \frac{\mathbbm{Pr}[{\chi^{\prime}_{1}}^2 (\lambda_i(\boldsymbol{Y}_{-i})) \ge  \bar{\Psi}_{1}^{-1}(\alpha_1)]}{\alpha_1} \mathbbm{1} \left (R(\boldsymbol{Y}_{-i}) = r \right ) \right \}, \nonumber
\end{eqnarray}
with the inequality following from Lemma 3.

Thus, for the BH method where $\alpha_r=r \alpha_1$, we finally have
\begin{eqnarray} \textrm {FDR} & \le & \sum_{i \in I_0} \mathbb{E} \left \{ \mathbbm{Pr}[{\chi^{\prime}_{1}}^2 (\lambda_i(\boldsymbol{Y}_{-i})) \ge  \bar{\Psi}_{1}^{-1}(\alpha_{1})] \sum_{r=0}^{d-1} \mathbbm{1} \left (R(\boldsymbol{Y}_{-i}) = r \right ) \right \} \nonumber \\
& = & \sum_{i \in I_0} \mathbb{E} \left \{ \mathbbm{Pr}[{\chi^{\prime}_{1}}^2 (\lambda_i(\boldsymbol{Y}_{-i})) \ge  \bar{\Psi}_{1}^{-1}(\alpha_{1}) ] \right \} = \sum_{i \in I_0} \mathbbm{Pr}[Y_i^2 \ge  \bar{\Psi}_{1}^{-1}(\alpha_{1})] \nonumber \\ & = & \sum_{i \in I_0} \mathbbm{Pr}[Z_i^2 \ge (1-R_i^2) \bar{\Psi}_{1}^{-1}(\alpha_{1})] \le \sum_{i=1}^d \bar{\Psi}_{1}\left (w_i\bar{\Psi}_{1}^{-1}(\alpha_{1})\right ) \nonumber \\ &  = & \alpha.
\end{eqnarray}
This proves the theorem.

\subsection{Wighted $t^2$ values} The setting for the multiple testing problem involving two-sided $t$-tests is as follows: Given  $\boldsymbol{X} = (X_1, \ldots, X_d)^{\prime} \sim {N}_d(\boldsymbol{\mu},\tau^2\boldsymbol{\Sigma})$, with
$\boldsymbol{\mu} = (\mu_1, \ldots, \mu_d)^{\prime}$ being an unknown vector and $\boldsymbol{\Sigma} = ((\sigma_{ij}))$ being a known positive definite matrix, and $V \sim \tau^2 \chi_m^2$, independently of $\boldsymbol{X}$, with $\tau^2 >0$ being an unknown scalar and $m >0$ being some known degrees of freedom, our problem is to test $H_i: \mu_i=0$ against $\mu_i \neq 0$, simultaneously for $i=1, \ldots, d$, subject to a control of the FDR at $\alpha$.

Let $\Psi_{m,n}$ denote the cdf of $F_{m,n}$, the random variable having central $F$ distribution with $m$ and $n$ degrees of freedom, and $\bar {\Psi}_{m,n} = 1 -{\Psi}_{m,n}$, Then $P_i = {\bar{\Psi}}^{-1}_{1,m}(T_i^2), \; \mbox{where} \; T_i= \frac{\sqrt{m}Z_{i}}{\sqrt{V}}, \; \mbox{for} \; i=1, \ldots, d$, are the $p$-values in the two-sided $t$-tests. As in the case of two-sided $z$-tests, we consider using the following increasing functions of these $p$-values:
\begin{eqnarray}\label{p-val2} \tilde{P}_i = \bar{\Psi}_{1,m} \left ( w_i^{-1}\bar {\Psi}_{1,m}^{-1} (P_i) \right ), \; \mbox{where} \; w_i = 1 -R_i^2, i=1, \ldots, d, \end{eqnarray} to run the BH method based on the critical constants $i\alpha_1$, $i=1, \ldots, d$, with an appropriately chosen $\alpha_1$ depending on the level of FDR control. In other words, we consider the weighted $t^2$ values $w_i^{-1}T_i^2 = \frac{mY_i^2}{V}$, $i=1, \ldots, d$, instead of the $t^2$ values $T_i^2$'s, in the BH method with the critical constants $\bar{\Psi}_{1,m}^{-1}((d-i+1)\alpha_1)$, $i=1, \ldots, d$ (as explained in Remark 1). This is our proposed weighted BH method in the two-sided $t$-tests. The fact that it controls the FDR is stated in the following theorem.

\vskip 6pt
\begin{theorem} \label{theorem2} The FDR of the BH method applied to the $\tilde{P}_i$'s in (\ref{p-val2}) using the critical constants $i\alpha_1$, $i=1, \ldots, d$, with $\alpha_1$ satisfying
\begin {eqnarray} \label{cond_theorem2}
& & \sum_{i=1}^d \bar{\Psi}_{1,m} \left ( w_i \bar {\Psi}_{1,m}^{-1} (\alpha_1) \right ) = \alpha, \end{eqnarray} controls the FDR at $\alpha$. \end{theorem}

The following two lemmas will paly a key role in our proof of this theorem.

\begin {lemma} \begin {eqnarray} \frac{\bar{\Psi}_n (\theta w)}{\bar{\Psi}_n (\theta w^{\prime})} \uparrow \theta > 0, \; \mbox{for any fixed} \; 0 < w < w^{\prime}< \infty. \end {eqnarray}
\end {lemma}

{\it Proof}. Let $W \sim \frac{1}{\theta} \chi_n^2$. Since the pdf of $W$ at $w$, which is $\theta^{\frac{n}{2}}e^{-\frac{1}{2}\theta w} w^{\frac{n}{2} -1}/2^{\frac{n}{2}}\Gamma(\frac{n}{2})$, is TP$_2$ (totally positive of order two) in $(w, \theta^{-1})$, the survival function of $W$, $\bar{\Psi}_n (\theta w) = \mathbb{P}(W \ge w)$, is also TP$_2$ in $(w, \theta^{-1})$. So, the required monotonicity result in the lemma holds. (See \cite{Karlin1968} for TP$_2$ related results).

\vskip 6pt

\begin{lemma} \begin {eqnarray} \frac{\bar{\Psi}_{m+h,n}({\bar{\Psi}}^{-1}_{m,n}(u))}
{u} \downarrow u\in (0,1), \nonumber \end{eqnarray}  for any fixed $m, n, h >0$. \end{lemma}

{\it Proof}. This lemma can be proved as in Lemma 3 by showing that the function
$g(u) = {\bar{\Psi}}_{m+h,n}({\bar{\Psi}}^{-1}_{m,n}(u)), \; u \in (0,1)$, for any fixed $m, n, h >0$,
is concave in $u\in(0,1)$.

\vskip 6pt

{\it Proof of Theorem 2}. We can assume without any loss of generality that $\tau^2 =1$ and that $\boldsymbol{\Sigma}$ is the correlation matrix. As in our proof of Theorem 1, we first see using the formula in (\ref{FDR1}) that the FDR of a step-up test based on the $p$-values in (\ref{p-val2}) and any set of critical constants $\alpha_i, i=1, \ldots, d$, is given by
\begin {eqnarray} \label{eqn1_theorem2} \textrm{FDR} & = & \sum_{i \in I_0} \sum_{r=0}^{d-1} \mathbb{E} \left \{ \frac{\mathbbm{Pr} [Y_i^2 \ge \frac{V}{m} \bar{\Psi}_{1,m}^{-1}(\alpha_{r+1})|~ \boldsymbol{Y}_{-i}, V ]\mathbbm{1} \left (R(\boldsymbol{Y}_{-i}, V) = r \right )}{r+1} \right \} \nonumber \\
& = & \sum_{i \in I_0} \sum_{r=0}^{d-1} \mathbb{E} \left \{ \frac{\mathbbm{Pr} [{\chi^{\prime}_1}^2 (\lambda_i(\boldsymbol{Y}_{-i}))\ge \frac{V}{m} \bar{\Psi}_{1,m}^{-1}(\alpha_{r+1})]}{r+1}~ \big |~ V)\mathbbm{1} \left (R(\boldsymbol{Y}_{-i}, V) = r \right ) \right \} \nonumber \\ & =  & \sum_{i \in I_0} \sum_{r=0}^{d-1} \mathbb{E} \left [ \mathbb{E}_{J_i} \left \{ \frac{\bar{\Psi}_{1+2J_i} \left ( \frac{V}{m} \bar{\Psi}_{1,m}^{-1}(\alpha_{r+1}) \right )}{r+1} ~\bigg |~ \boldsymbol{Y}_{-i}, V \right \} \mathbbm{1} \left (R(\boldsymbol{Y}_{-i}, V) = r \right )\right ]. \end{eqnarray}

The right-hand side of (\ref{eqn1_theorem2}) can be re-written, as in (\ref{FDR2}), to have the following:

\begin {eqnarray} & & \textrm{FDR} - \sum_{i \in I_0}\bar{\Psi}_{1,m} ((1-R_i^2)\bar{\Psi}_{1,m}^{-1}(\alpha_1))\nonumber \\
& \le & \sum_{i \in I_0} \sum_{r=1}^{d-1} \mathbb{E} \left [ \mathbb{E}_{J_i} \left \{ \frac{\bar{\Psi}_{1+2J_i} (\frac{V}{m}\bar{\Psi}_{1,m}^{-1}(\alpha_{r+1}))}{r+1}- \frac{\bar{\Psi}_{1+2J_i} (\frac{V}{m}\bar{\Psi}_{1,m}^{-1}(\alpha_{r}))}{r}~\bigg |~\boldsymbol{Y}_{-i}, V \right \} \times \right. \nonumber \\
& & \qquad \qquad \left.  \mathbbm{1} \left (R(\boldsymbol{Y}_{-i}, V) \ge r \right ) \right ] \nonumber \\
& = & \sum_{i \in I_0} \sum_{r=1}^{d-1} \mathbb{E} \left ( \mathbb{E}_{V} \left [ \left \{ \frac{\bar{\Psi}_{1+2J_i} (\frac{V}{m}\bar{\Psi}_{1,m}^{-1}(\alpha_{r+1}))}{(r+1)\bar{\Psi}_{1+2J_i} (\frac{V}{m}\bar{\Psi}_{1,m}^{-1}(\alpha_{r}))} - \frac{1}{r} \right \} \times  \right. \right. \nonumber \\ & &\qquad \qquad  \left. \left. \mathbbm{1} \left (R(\boldsymbol{Y}_{-i}, V) \ge r \right )\bar{\Psi}_{1+2J_i} \left (\frac{V}{m}\bar{\Psi}_{1,m}^{-1}(\alpha_{r}) \right ) ~|~\boldsymbol{Y}_{-i}, J_i \right ] \right ). \nonumber \end {eqnarray}

From Lemma 4, we note that, given ($\boldsymbol{Y}_{-i}, J_i)$, $$\frac{\bar{\Psi}_{1+2J_i} (\frac{V}{m}\bar{\Psi}_{1,m}^{-1}(\alpha_{r+1}))}{(r+1)\bar{\Psi}_{1+2J_i} (\frac{V}{m}\bar{\Psi}_{1,m}^{-1}(\alpha_{r}))} - \frac{1}{r} \uparrow V >0,$$ Also, $\mathbbm{1} \left (R(\boldsymbol{Y}_{-i}, V) \ge r \right ) \downarrow V$, given ($\boldsymbol{Y}_{-i}, J_i)$. Thus, from Kimball's inequality, we see that
\begin {eqnarray}
& &  \mathbb{E}_{V} \left [ \left \{ \frac{\bar{\Psi}_{1+2J_i} (\frac{V}{m}\bar{\Psi}_{1,m}^{-1}(\alpha_{r+1}))}{(r+1)\bar{\Psi}_{1+2J_i} (\frac{V}{m}\bar{\Psi}_{1,m}^{-1}(\alpha_{r}))} - \frac{1}{r} \right \} \times  \right. \nonumber \\
& & \qquad \qquad \left. \mathbbm{1} \left (R(\boldsymbol{Y}_{-i}, V) \ge r \right )\bar{\Psi}_{1+2J_i} (\frac{V}{m}\bar{\Psi}_{1,m}^{-1}(\alpha_{r}))~\bigg |~\boldsymbol{Y}_{-i},J_i \right ] \nonumber \\
& \le & \left \{ \frac{ \mathbb{E}_{V} \left [ \bar{\Psi}_{1+2J_i} (\frac{V}{m}\bar{\Psi}_{1,m}^{-1}(\alpha_{r+1}))\right ] }{r+1} - \frac{ \mathbb{E}_{V} \left [\bar{\Psi}_{1+2J_i} (\frac{V}{m}\bar{\Psi}_{1,m}^{-1}(\alpha_{r})) \right ]}{r} \right \} \times  \nonumber \\
& & \qquad \qquad \frac{ \mathbb{E}_{V} \left [ \mathbbm{1} \left (R(\boldsymbol{Y}_{-i}, V) \ge r \right )\bar{\Psi}_{1+2J_i} (\frac{V}{m}\bar{\Psi}_{1,m}^{-1}(\alpha_{r})) \bigg |~\boldsymbol{Y}_{-i},J_i \right ]}{\mathbb{E}_{V} \left [\bar{\Psi}_{1+2J_i} (\frac{V}{m}\bar{\Psi}_{1,m}^{-1}(\alpha_{r})) \bigg |~\boldsymbol{Y}_{-i},J_i \right ] }. \nonumber \end{eqnarray}
For the BH method, where $\alpha_r = r\alpha_1$, this is less than or equal to zero, for each $r=1, \ldots, d$, since
\begin {eqnarray}
& & \frac{ \mathbb{E}_{V} \left [ \bar{\Psi}_{1+2J_i} (\frac{V}{m}\bar{\Psi}_{1,m}^{-1}(\alpha_{r+1})) \right] }{r+1} - \frac{ \mathbb{E}_{V} \left [\bar{\Psi}_{1+2J_i} (\frac{V}{m}\bar{\Psi}_{1,m}^{-1}(\alpha_{r})) \right ]}{r} \nonumber \\
&= & \alpha_1 \left [ \frac{ \bar{\Psi}_{1+2J_i,m} (\bar{\Psi}_{1,m}^{-1}(\alpha_{r+1}))}{\alpha_{r+1}} -
\frac{ \bar{\Psi}_{1+2J_i,m} (\bar{\Psi}_{1,m}^{-1}(\alpha_{r}))}{\alpha_r} \right ] \nonumber \\
& \le & 0, \nonumber \end{eqnarray} from Lemma 5. This proves the required inequality for the BH method: $$ \textrm{FDR} \le \sum_{i \in I_0}\bar{\Psi}_{1,m} (w_i\bar{\Psi}_{1,m}^{-1}(\alpha_1)) \le \alpha $$ by considering $\alpha_1$ satisfying (\ref{cond_theorem2}). This completes the proof of Theorem 2.

\section{Further discussion}

Here, we discuss the aforementioned main results in some special cases.

Clearly, when the covariance matrix of the underlying multivariate normal distribution is diagonal with known or unknown scalar multiple of known entries, our proposed weighted BH methods reduce to the corresponding original BH methods involving the usual, un-weighted $z^2$- or $t^2$-values. For other cases, such as equi-correlated multivariate normal and the multivariate normal arising in variable selection under linear regression, the proposed methods are presented below.

\subsection{Equi-correlated multivariate normal}
Let the $\boldsymbol{X}$ in Sections 3.1 and 3.2 have the following correlation matrix: $(1-\rho)I_d+\rho \boldsymbol{1}_d\boldsymbol{1}_d^{\prime}$ [with $\boldsymbol{1}_d = (1, \ldots, 1)^{\prime}: d \times 1$], for some $-1/(p-1) < \rho < 1$. In this case, all $d$ multiple correlations are same as $R^2 = (p-1)\rho^2/[1+(p-2)\rho]$, and so the weighted BH methods in Theorems 1 and 2 are based on the corresponding $\tilde{P}_i$'s with
$$w_i = \frac{(1-\rho)[1+(p-1)\rho]}{1+(p-2)\rho}, \; \forall \; i=1, \ldots, d.$$

\subsection{FDR controlled variable selection}
Consider the variable/feature selection problem under the following linear regression model:
\begin {eqnarray}\label{VS_Model} \boldsymbol{Y} = \boldsymbol{X}\boldsymbol{\beta} + \boldsymbol{\epsilon}, \end {eqnarray} where $\boldsymbol{Y}$ is $n$-dimensional response vector, $\boldsymbol{X} = (\boldsymbol{X}_1, \ldots, \boldsymbol{X}_d)$ is $n \times d$ design matrix of rank $d \le n$ with its columns representing the known vectors of observations on the $d$ variables/features $X_1, \ldots, X_d$, $\boldsymbol{\beta} = (\beta_1, \ldots, \beta_d)$ is the unknown vector of regression coefficients corresponding to these variables/features, and $\boldsymbol{\epsilon} \sim N_d(\boldsymbol{0}, \tau^2 \boldsymbol{I}_d)$ is the Gaussian noise.

Variable selection can be framed as a multiple testing problem, where the null hypothesis $H_i: \beta_i=0$ is tested against its  alternative $\beta_i \neq 0$, simultaneously for $i=1, \ldots, d$, and the variables corresponding to the rejected nulls, according to a multiple testing procedure based on some estimates of the regressing coefficients, are selected/discovered as the important variables. The FDR of the multiple testing procedure would be a powerful measure of potential errors in the selection.

The ordinary least squares estimate of ${\boldsymbol{\beta}}$,  given by $\hat{\boldsymbol{\beta}}= (\hat{\beta}_1, \ldots, \hat{\beta}_d)^{\prime} =  \boldsymbol{A}^{-1}\boldsymbol{X}^{\prime}\boldsymbol{Y}$, where $\boldsymbol{A} = \boldsymbol{X}^{\prime}\boldsymbol{X}$, is distributed as $N_d(\boldsymbol{\beta}, \tau^2 \boldsymbol{A}^{-1})$, and is independent of ${\hat{\tau}}^2 = [\|\boldsymbol{Y}\|^2 - \hat{\boldsymbol{\beta}}^{\prime} \boldsymbol{A} \hat{\boldsymbol{\beta}}]/(n-d) \sim \tau^2 \chi^2_{n-d}/(n-d)$.
Therefore, a natural choice for a powerful FDR controlling procedure that can be used for variable selection would be the BH method based on the following $t^2$-values: \begin {eqnarray} T_i^2 = \frac{Z_i^2}{\hat{\tau}^2}, \; \mbox{where} \; Z_i^2 = \frac{\hat{\beta}_i^2}{a^{ii}}, \; \mbox{for} \; i=1, \ldots, d, \end {eqnarray}
with $a^{ii}$ being the $i$th diagonal entry of $\boldsymbol{A}^{-1}$. Unfortunately, as noted in Introduction, there is no theoretical guarantee that the BH method (in its original form) can control the FDR.

Following the development of our proposed method in Section 3.2, we can now propose a novel BH-type step up method for variable selection with proven FDR control. To that end, we first note that $R_i^2$, the squared multiple correlation between $\hat{\beta}_i$ and $\hat{\boldsymbol{\beta}}_{-i}$ equals $1- a_{ii}^{-1}/a^{ii}$. So, so we can formally describe our proposed FDR controlling procedure in variable selection as follows:

\begin {definition} [Weighted BH method for variable selection] Let $P_i= \bar{\Psi}_{1,n-d}(T_i^2)$ be the original p-values corresponding to $H_i$. Run the BH-type stepup method using $\tilde{P}_i=  \bar{\Psi}_{1,n-d}(w_i^{-1}\bar{\Psi}^{-1}_{1,n-d}(P_i))$, where $w_i=a_{ii}^{-1}/a^{ii}$, $i=1, \ldots, d$, and the critical constants $i\alpha_i$, $i=1, \ldots, d$, with $\alpha_1$ being such that
$$\sum_{i=1}^d \bar{\Psi}_{1,n-d}(w_i\bar{\Psi}^{-1}_{1,n-d}(\alpha_1)) = \alpha.  $$ \end {definition}

\section{Concluding remarks} This paper answers the following question, a paraphrase of what we have stated in the abstract: Can the BH method control, or a BH-type stepup method be developed to control, the FDR in the two-sided $z$- or $t$-tests for multiple testing of multivariate normal means? The proposed weighted BH methods in Section 3 provide an affirmative answer to this question, assuming of course that the correlations are known. When the correlations, as well as the variances, are unknown, answering this question for two-sided $t$-tests based on the marginal $t^2$-statistics in the Hotelling's T$^2$ test is an important, yet challenging, open problem.

The proposed weighted BH methods can be viewed as some sorts of adjustment of the BH method to the underlying correlation structure, like the dependence adjusted BH (DBH) method in \cite{Fithian2022}. However, it is important to point out that, while the DBH fully captures the underlying correlation structure and is quite powerful, it could be less user-friendly in practical applications, as it is not expressible in a closed form. It is implemented through a computer assisted algorithm requiring extensive computation. Regarding the novelty of our proposed weighted BH method involving two-sided $t$-tests in FDR controlled variable selection (in Section 4.2), we must point out that, while the knockoff based FDR controlling methods in \cite{Barber2015} and \cite{Sarkar2022} offer powerful alternatives, they are not applicable when $d \le n < 2d$; they require $n$ to be greater than or equal to 2$d$.



\end{document}